\documentclass[12pt]{amsart}
\usepackage{amssymb,amsmath,amsthm}
\newtheorem{lemma}{Lemma}

\newtheorem{thm}{Theorem}

\def\d{\displaystyle}
\begin{document}
\title
{Weak type radial convolution operators on free group}
\author{\frame{\ Tadeusz Pytlik } and Ryszard Szwarc}
\address{Institute of Mathematics,
University of Wroc\l aw, pl.\ Grunwaldzki 2/4, 50-384 Wroc\l aw, Poland \newline  and
\newline \indent Institute of Mathematics and Computer Science, University of Opole, ul. Oleska 48,
45-052 Opole, Poland}
 \email{szwarc2@gmail.com}
\date{}
\thanks{The second author was supported by European Commission Marie Curie Host
Fellowship for the Transfer of Knowledge ``Harmonic Analysis, Nonlinear Analysis and
Probability''  MTKD-CT-2004-013389 and by MNiSW Grant N201 054 32/4285}
\begin{abstract}
Radial convolution operators on free groups with nonnegative kernel of weak type $(2,2)$
and of restricted weak type $(2,2)$ are characterized. Estimates of weak type $(p,p)$ are
obtained as well for $1<p<2.$
\end{abstract}

\maketitle
\section{Introduction}
A discrete group $G$ is called amenable if there exists a linear functional $m$
on $\ell_{\mathbb{R}}^{\infty}(G)$ such that
\begin{enumerate}
\item[(1)] $\d\inf_{x\in G}f(x)\le m(f)\le \sup_{x\in G}f(x),$
\item[(2)] $m(_xf)=m(f), \quad{\rm where}\  _xf(y)=f(x^{-1}y).$
\end{enumerate}
$m$ is called a left invariant mean. Then the functional $M(f)=m(m(f_x))$ satisfies (1), (2) and is
also right invariant, where $f_x(y)=f(yx).$

Let $G$ be a discrete group. Consider a symmetric probability measure
$\mu$ on $G,$ i.e.
\begin{align*}
&\mu=\sum_{x\in G} \mu(x)\delta_x, & &\mu(x)
\ge 0,\\
& \sum_{x\in G}\mu(x)=1, && \mu(x^{-1})=\mu(x).
\end{align*}
The left convolution operator $\lambda(\mu)$ with $\mu$ is
bounded  on $\ell^2(G)$ and
$$\|\lambda(\mu)(f)\|_2=\| \mu*f\|_2\le \|f\|_2,\quad f\in  \ell^2(G).$$
Indeed
$$\| \mu*f\|_2=\left\|\sum_{x\in G}\mu(x)[ \delta_x*f]\right\|_2
\le\sum_{x\in G} \mu(x)\| \delta_x*f\|_2=\|f\|_2.
$$
Thus $\|\lambda(\mu)\|_{2\to 2}\le 1.$

 Kesten \cite{k1} showed that a discrete group $G$ is amenable iff for any symmetric
probability measure $\mu$ on $G$ we have $\|\lambda(\mu)\|_{2\to 2}=1. $
 He showed that $G$ is amenable if
condition is satisfied for one measure $\mu$ such that ${\rm supp}\,\mu$
generates $G$ algebraically. In particular let $G$ be generated by $g_1,g_2,\ldots, g_k$
and  $\mu={1\over 2k}\d\sum_{i=1}^k(\delta_{g_i}+\delta_{g_i^{-1}}).$
Then $G$ is amenable if and only if $\|\lambda(\mu)\|_{2\to 2}=1.$

In \cite{f} F\o lner came up with  another property  equivalent to amenability.
We say that a discrete group $G$ satisfies the F\o lner condition if for any number $\varepsilon>0$ and any finite set $K\subset G$ there
exists a finite set $N\subset G$ such that
\begin{equation}\label{fo}
|xN\bigtriangleup N|< \varepsilon |N|,\qquad x\in K.
\end{equation}
In other words $N$ is almost $K$ invariant.
He showed that
 $G$ is amenable if and only if the F\o lner condition holds.

Assume that $G$ is amenable. Let $\mu$ be a probability measure  with finite
support $K.$ For $\varepsilon=\eta^2 >0$ choose $N$ so as to satisfy (\ref{fo}). Then
\begin{multline*}
\|\mu*\chi_N-\chi_N\|_2= \left\|\sum_{x\in K}\mu(x)[\chi_{xN}-\chi_N]\right
\|_2
\le \sum_{x\in K}\mu(x)\|\chi_{xN}-\chi_N\|_2\\=\sum_{x\in
K}\mu(x)\|\chi_{xN\bigtriangleup N }\|_2
=\sum_{x\in K}\mu(x)|xN\bigtriangleup N|^{1/2}\le
\eta|N|^{1/2}=\eta\|\chi_N\|_2.
\end{multline*}
Therefore
\begin{multline*}\langle \mu*\chi_N,\chi_N\rangle_{\ell^2(G)}\\=
\langle  \chi_N,\chi_N\rangle_{\ell^2(G)}+\langle \mu*\chi_N-\chi_N,\chi_N\rangle_{\ell^2(G)}
\ge
 (1-\eta)\|\chi_N\|^2_2,
 \end{multline*}
which implies
\begin{equation}\label{fo2}
\sup_{N,M-{\rm finite}}{\langle\mu*\chi_N,\,\chi_M\rangle\over
\|\chi_N\|_2\|\chi_M\|_2}=1=\|\lambda(\mu)\|_{2\to 2}.
\end{equation}
The same holds (with the same proof) for any $1< p<+\infty,$
 i.e.
\begin{equation}\label{fop}\sup_{N,M-{\rm finite}}{\langle\mu*\chi_N,\, \chi_M\rangle \over
\|\chi_N\|_p\|\chi_M\|_{p'}}=1=\|\lambda(\mu)\|_{p\to p},
\end{equation}
where $p'=p/(p-1).$

We will use the notion of $L^{p,q}$ spaces, which have been introduced by Lorentz (see \cite{bl}).
Consider a  general $\sigma$-finite measure space $(\Omega,\omega)$ and $1< p<+\infty.$
For $f\in L^p(\Omega,\omega)$ and $t>0$ we have
$$t^p\,\omega\{x:|f(x)|>t\}\le \int_{\Omega}|f(x)|^pd\omega(x).$$
  Functions for which the left hand side is bounded form a linear space
   $$L^{p,\infty}(\Omega,\omega)=
\left \{f:\sup_{t>0}t^p\omega\{x:|f(x)|>t\}<+\infty\right\}.$$
called the weak $L^p$ space. This space contains $L^p(\Omega,\omega).$

For $p'=p/(p-1)$ the predual of $L^{p',\infty}(\Omega,\omega)$ with respect to the standard inner product
is denoted by $L^{p,1}(\Omega,\omega).$
We have
   $$L^{p,1}(\Omega,\omega)\subset
 L^p(\Omega,\omega)\subset L^{p,\infty}(\Omega,\omega).$$
For $p>1$ these spaces are normed.

 Any linear operator mapping $L^p$ into itself is called of
 strong type $(p,p).$ Linear operators $T$ mapping $L^{p}(\Omega,\omega)$ into\
 $L^{p,\infty}(\Omega,\omega)$ are
 called of weak type $(p.p),$ while those which map
 $L^{p,1}(\Omega,\omega)$ into  $L^{p,\infty}(\Omega,\omega)$
 are called of restricted weak type $(p,p).$

   We will use the following facts.  A linear operator $T$
 is bounded from $L^{p,1}$ into a Banach space $X$ if and only if
\begin{equation}\label{weak}
\|T\|_{L^{(p,1)}\to X}=\sup_{E\subset \Omega}{\|T\chi_E\|_X\over \|\chi_E\|_p}<+\infty.
\end{equation}
A linear operator
$T$
 is bounded from $L^{p,1}$ into $L^{p,\infty}$ if and only if
\begin{equation}\label{rweak}
\|T\|_{(p,1)\to (p,\infty)}=\sup_{E, F\subset \Omega}{|\langle T\chi_E,\,\chi_F\rangle|\over \|\chi_E\|_p
\|\chi_F\|_{p'}}<+\infty.
\end{equation}
Using this and duality between spaces $L^{(p',1)}$ and $L^{(p,\infty)}$ we obtain
\begin{equation}\label{dual}
\|T\|_{p\to (p,\infty)}=\|T^*\|_{(p',1)\to p'}=
\sup_{E\subset \Omega}{\|T^*\chi_E\|_{p'}\over \|\chi_E\|_{p'}}.
\end{equation}
The equalities (\ref{fo2}) and (\ref{fop}) can be interpreted as follows.
If  the group $G$ is discrete and amenable and $\mu$ is a symmetric probability measure
on $G,$ then
\begin{equation}\label{duality}
\|\lambda(\mu)\|_{p\to p}= \|\lambda(\mu)\|_{(p',1)\to p'}=\lambda(\mu)\|_{p\to (p,\infty)}=\|\lambda(\mu)\|_{(p,1)\to (p,\infty)}=1.
\end{equation}
Hence for these groups convolution operators with   nonnegative
functions  of strong type $(p,p),$  of weak type $(p,p)$ and
of restricted weak type $(p,p)$ coincide for any $1<p<\infty.$

The situation is entirely different for nonamenable groups.
Only special examples have been studied. It has been  shown \cite{sz} that for
$p=2$ and $G= \mathbb{F}_k,$ the free group on $k$ generators, $k\ge 2,$ there exist nonnegative functions
$f$
on $G$ such that $\|\lambda(f)\|_{2\to (2,\infty)}$ is finite while $\|\lambda(f)\|_{2\to 2}$ is infinite,
i.e. there exist convolution operators with nonnegative functions of weak type (2,2) which are not
of strong type (2,2).
The same has been shown for $1<p<2$ \cite{sz2}. These functions $f$ can be chosen to be radial, i.e.
constant on elements of the group $G$ of the same  length. It is an open problem if these results remain
true for any discrete nonamenable group.

In this work will focus on $G=\mathbb{F}_k.$ We are going to determine all nonnegative
radial functions $f$ on $G$ such that $\lambda(f)$ is of weak type (2,2), as well those
$f$ for which $\lambda(f)$ is of restricted weak type (2,2). In particular we prove that
these spaces are different. Next we will turn our attention to the case $1<p<2.$ By using
interpolation machinery, duality and the results for $p=2$ we will be able to determine
nonnegative radial functions $f$ for which $\lambda(f)$ is of weak type $(p,p).$ In this
way we obtain a simpler proof of the upper estimate of $\|\lambda(f)\|_{p\to (p,\infty)}$
obtained in \cite{cms}. Our method   does not rely on deep theorems on representation
theory
\section{Radial convolution operators of weak type $(2,2)$}
Let $\mathbb{F}_k={\rm gp}\{g_1,g_2,\ldots, g_k\}$ be a free group on $k\ge 2$
generators. The  group consists of reduced words in generators and their inverses. This
representation is unique. The number of letters in reduced form defines length function
on $ \mathbb{F}_k.$ Let $\chi_n$ denote the indicator function of words of length $n.$
There are $2k(2k-1)^{n-1}$ such words. as we have $2k$ choices for the first letter and
$2k-1$ choices for every consecutive one. Let $q=2k-1.$ The next theorem generalizes the
estimate for $\|\lambda(\chi_n)\|_{2\to (2,\infty)}$ given in \cite{sz}.

\begin{thm}\label{weak1}
Let $f=\sum_{n=0}^\infty f_n\chi_n.$  The operator
$\lambda(f)$ is of weak type $(2,2)$   if
$$A(f):=\sum_{n,m=0}^\infty |f_n||f_m|q^{-(n+m)/2}\{1+\min(n,m)\}<+\infty.$$
Moreover if $f_n\ge 0$ the condition is necessary and
$${1\over 6}A(f)\le \|\lambda(f)\|^2_{2\to (2,\infty)}
\le 4A(f).$$
\end{thm}

{\it Proof.}
By (\ref{duality}), instead of estimating $\|\lambda(f)\|_{2\to (2,\infty)}$ we may estimate
$\|\lambda(f)\|_{(2,1)\to 2},$ which (see (\ref{weak})) is equivalent to
$$\sup_{E\subset \mathbb{F}_r}{\|f*\chi_E\|_2\over |E|^{1/2}}.$$
We have
$$
\|f*\chi_E\|_2^2=\langle f*f*\chi_E*,\chi_E\rangle=
\sum_{n,m=0}^\infty f_nf_m\langle \chi_n*\chi_m*\chi_E, \chi_E\rangle.
$$
Simple calculation gives that for $n\ge 1$ we have
$$
\chi_n*\chi_m=q^{n-1}\delta_{n}^m\chi_0+
\sum^{n+m}_{k=|n-m|\atop k\equiv n+m \mod 2}q^{(n+m-k)/2}\,\chi_k
.
$$ Clearly $\chi_0*\chi_0=\chi_0.$
Therefore
$$\chi_n*\chi_m\le   2\sum^{n+m}_{k=|n-m|\atop k\equiv n+m \mod 2}q^{(n+m-k)/2}\chi_k.$$
Hence
$$
\|f*\chi_E\|_2^2\le 2\sum_{n,m=0}^\infty f_nf_m\,
q^{(n+m)/2}\sum^{n+m}_{k=|n-m|\atop k\equiv n+m \mod 2}q^{-{k/
2}}\,
\langle \chi_k*\chi_E,\chi_E\rangle.
$$
\begin{lemma}\label{w22} $$\langle \chi_{k}*\chi_E,\chi_E\rangle \le
2q^{[k/2]}|E|.$$
\end{lemma}

{\it Proof.}
Define an operator $P_k$ by the rule
$$\langle P_k\delta_x,\delta_y\rangle=
\begin{cases}\langle \chi_{k}*\delta_x,\delta_y\rangle & {\rm if}\  |x|\ge |y|\\
0 & {\rm if}\  |x|<|y|.
\end{cases}$$
Then $$\langle \chi_{k}*\delta_x,\delta_y\rangle\le \langle
P_k\delta_x,\delta_y\rangle+\langle
\delta_x,P_k\delta_y\rangle.$$
This implies
$$\langle \chi_{k}*\chi_E,\chi_E\rangle \le
2\langle P_k\chi_E,\chi_E\rangle\le 2\|P_k\chi_E\|_1\le 2|E|\,\sup_x\|P_k\delta_x\|_1
 $$
 Next $$P_k\delta_x=\sum_{ |w|=k\atop |wx|\le|x|}\delta_{wx}.$$
Let $w=w_1w_2$ where $|w_1|\le |w_2|\le (k+1)/2.$ The conditions $|w|=k$ and $|wx|\le
|x|$ imply that $w_2$ is determined by the first $[(k+1)/2]$ letters of $x.$
Hence we have as many terms in the sum as choices for $w_1,$ i.e. at most
  $q^{[k/2]}.$ Thus
$$\|P_k\delta_x\|_1\le q^{[k/2]}.$$
Therefore
$$\langle \chi_{k}*\chi_E,\chi_E\rangle \le 2 q^{[k/2]}|E|.$$
\qed

\noindent Lemma 1 implies that
\begin{multline*}
{\|f*\chi_E\|_2^2\over |E|}\le 4\sum_{n,m=0}^\infty |f_n||f_m| q^{(n+m)/2}
\sum^{n+m}_{k=|n-m|\atop k\equiv n+m \mod 2}\\=
4 \sum_{n,m=0}^\infty |f_n||f_m| q^{(n+m)/2} \{1+\min(m,n)\}.
\end{multline*}
We obtain the upper estimate
$$\|\lambda(f)\|^2_{2\to (2,\infty)}\le 4 \sum_{n,m=0}^\infty |f_n||f_m| q^{(n+m)/2}
\{1+\min(m,n)\}.$$
On the other hand if $f_n\ge 0$ we have
\begin{multline*}
\|\lambda(f)\|_{2\to (2,\infty)}^2\ge
{q\over q+1}q^{-2k}\|f*\chi_{2k}\|^2_2\ge {2\over 3}\,q^{-{2k}}\left\|\sum_{n=0}^\infty
f_n(\chi_n*\chi_{2k})\right \|_2^2\\
\ge {2\over 3}\,q^{-{2k}}\left\|\sum_{n=0}^\infty
f_n \sum^{n+2k}_{l=|n-2k|\atop l\equiv n \mod 2}q^{(n+2k-l)/2}\,\chi_l\right
\|_2^2\\={2\over 3}\,\left\|\sum_{l=0}^\infty
q^{-l/2}\chi_l\left (
 \sum^{2k+l}_{n=|2k-l|\atop n\equiv  l \mod 2}f_nq^{n/2}\right )\right\|_2^2
 \ge {2\over 3}\,\sum_{l=0}^\infty \left (
 \sum^{2k+l}_{n=|2k-l|\atop n\equiv  l \mod 2}f_nq^{n/2}\right )^2\\
 \ge {2\over 3}\,\sum_{l=0}^{2k} \left (
 \sum^{2k+l}_{n=2k-l\atop n\equiv  l \mod 2}f_nq^{n/2}\right )^2\ge {2\over
 3}\,
 \sum_{n,m=0}^{2k} f_nf_m q^{(n+m)/2}\sum_{l=\max(2k-n,2k-m)\atop l\equiv n\equiv m  \mod 2}^{2k}.
\end{multline*}
 Considering even or odd values of $m$ and $n$ gives
\begin{eqnarray*}
\|\lambda(f)\|_{2\to (2,\infty)}^2&\ge & {2\over 3}\,\sum_{n,m=0}^k f_{2n}f_{2m}q^{n+m}
\{1+\min(n,m)\},\\
\|\lambda(f)\|_{2\to (2,\infty)}^2&\ge &{2\over 3}\, \sum_{n,m=0}^{k-1} f_{2n+1}f_{2m+1}q^{n+m+1}
\{1+\min(n,m)\}.
\end{eqnarray*}
Since $k$ is arbitrary
$$
\|\lambda(f)\|_{2\to (2,\infty)}^2 \ge {1\over 3}
\sum_{n,m=0\atop n\equiv m\mod 2}^\infty f_nf_m q^{(n+m)/2}\{1+\min(n,m)\}
 .$$
This implies
$$\|\lambda(f)\|_{2\to (2,\infty)}^2 \ge {1\over 6}\sum_{n,m=0 }^\infty f_nf_m
q^{(n+m)/2}\{1+\min(n,m)\},
$$
because the matrix $a(n,m)=1+\min (n,m)$ is positive definite. \qed
\begin{thm}
For $n\ge 0$ there holds
$$\|\lambda(\chi_n\|_{(2,1)\to (2,\infty)}\le cq^{n/2}.$$
\end{thm}
{\it Proof.}
We have
$$ \|\lambda(\chi_n\|_{(2,1)\to (2,\infty)}=\sup_{E,F\subset \mathbb{F}_r}
{\langle \chi_n*\chi_E,\chi_F\rangle\over |E|^{1/2}|F|^{1/2}}.$$
The proof will be completed if we show
\begin{equation}
\label{r22}
\langle \chi_n*\chi_E,\chi_F\rangle\le cq^{n/2}|E|^{1/2}|F|^{1/2}.
\end{equation}
We will prove (\ref{r22}) by modification of the argument used in the
proof
of Lemma \ref{w22}. Fix $\alpha\in \mathbb{R}.$ Let $Q_n^\alpha$ denote the operator defined by the rule
$$\langle Q_n^\alpha\delta_x,\delta_y\rangle=
\begin{cases}\langle \chi_{n}*\delta_x,\delta_y\rangle & {\rm if}\  |x|\ge q^\alpha|y|\\
0 & {\rm if}\  |x|<q^\alpha|y|.
\end{cases}$$
Then $$\langle \chi_{n}*\delta_x,\delta_y\rangle\le \langle
Q_n^\alpha\delta_x,\delta_y\rangle+\langle
\delta_x,Q_n^{-\alpha}\delta_y\rangle.$$
This implies
\begin{eqnarray}\label{basic}\langle \chi_{n}*\chi_E,\chi_F\rangle &\le&
  \|Q_n^\alpha\chi_E\|_1+\|Q_n^{-\alpha}\chi_F\|_1 \nonumber\\
&\le &|E|\, \sup_x\|Q_n^\alpha \delta_x\|_1 +|F|\, \sup_x\|Q_n^{-\alpha} \delta_x\|_1
  \end{eqnarray}
Next $$Q_n^\alpha\delta_x=\sum_{ |w|=n\atop |wx|\le q^{-\alpha}|x|}\delta_{wx}.$$
Let $w=w_2w_1$ where $|w_1|=[n/2]+[\alpha]$ and
$ |w_2|= n-[n/2]-[\alpha].$ The conditions $|w|=n$ and $|wx|\le
q^{-\alpha}|x|$ imply that $w_1$ is determined by the first $ [n/2]+[\alpha] $ letters of $x.$
Hence we have as many terms in the sum as choices for $w_2,$ i.e. at most
  $ q^{n-[n/2]-[\alpha]}.$ Thus
\begin{equation}\label{estim1} \|Q_n^\alpha\delta_x\|_1\le q^{3/2} q^{-\alpha}q^{n/2}.
\end{equation}
Similarly
$$\|Q_n^{-\alpha}\delta_x\|_1\le q^{3/2}q^{ \alpha}q^{n/2}.
$$
Hence by (\ref{basic}) we get
$$\langle \chi_{n}*\chi_E,\chi_F\rangle \le q^{3/2}q^{n/2}\{q^{-\alpha} \,|E|+
  q^{ \alpha} \,|F|\}. $$
Choosing $\alpha=(\log|E|-\log|F|)/(2\log q)$ gives
$$\langle \chi_{n}*\chi_E,\chi_F\rangle \le
2q^{3/2}q^{n/2}|E|^{1/2}|F|^{1/2}.$$
\qed
\begin{thm}\label{restr}
Let $f=\sum_{n=0}^\infty f_n\chi_n$ and $f_n\ge 0.$ The operator
$\lambda(f)$ is of restricted weak type $(2,2)$ if and only if
$f\in L^{(2,1)}.$
\end{thm}
{\it Proof.} By Theorem 2 we have
$$\|\lambda(\chi_n)\|_{(2,1)\to (2,\infty)}\le Cq^{n/2}$$
for some constant $C>0.$ Let $f=\sum_{n=0}^\infty f_n\chi_n.$
Then   triangle inequality yields
$$\|f\|_{(2,1)\to (2,\infty)}\le C\sum_{n=0}^\infty f_nq^{n/2}.$$
By \cite[Lemma 1]{p}
\begin{equation}\label{equiv}\sum_{n=0}^\infty f_nq^{n/2}\approx \|f\|_{(2,1)}.
\end{equation}
 On the other hand for $f_n\ge 0$ we have
\begin{multline*}
\|f\|_{(2,1)\to (2,\infty)}\ge C\sup_{n,m}\,q^{-(n+m)/2}\,\langle
f*\chi_n,\chi_m\rangle\\=C\sup_{n,m}\,q^{-(n+m)/2}\,\langle
f,\chi_m*\chi_n\rangle
\ge C\sum^{n+m}_{k=|n-m|\atop k\equiv n+m\mod 2}q^{k/2}f_k.
\end{multline*}
Taking $m=n$ or $m=n+1$ and letting $n$ tend to infinity gives
\begin{eqnarray*}
\|f\|_{(2,1)\to (2,\infty)}&\ge & C \sum_{k=0}^\infty q^{2k/2}f_{2k},\\
\|f\|_{(2,1)\to (2,\infty)}&\ge & C \sum_{k=0}^\infty
q^{(2k+1)/2}f_{2k+1}.
\end{eqnarray*}
Therefore $ \sum_{k=0}^\infty q^{k/2}f_k <+\infty ,$ i.e. $f\in
L^{(2,1)}$ by (\ref{equiv}).\qed
\section{Weak type $(p,p)$ for $1<p<2$}
Part of the next theorem, namely the first inequality is known from \cite{cms}.
Actually it has been simply observed there that the inequality follows by applying
 multilinear interpolation theorem to Pytlik's estimate for $\|\sum f_n\lambda(\chi_n)\|_{p\to p}$
 given in
\cite{p}. We will reprove the second inequality by applying the same
interpolation theorem to restricted weak type estimates given in the previous
section. In this way we skip $p\to p$ estimates
whose proof as given in \cite{p} is tricky, and later proof
given in \cite{cms} makes use of advanced representation theory.

\begin{thm} For $1<p<2$ and for $f=\sum_{n=0}^\infty f_n\chi_n$ we have
$$\|\lambda(f)\|_{p\to (p,\infty)}
 \le C\|f\|_{(p,p')}.$$
 Moreover if $f\ge 0$ then
 $$c\|f\|_{(p,p')}\le \|\lambda(f)\|_{p\to (p,\infty)}.$$
\end{thm}

{\it Proof.} The subscript $r$ will denote the subspace of radial functions, i.e.
functions of the form $\sum_{n=0}^\infty f_n \chi_n,$ where $f_n$ are complex coefficients.
By Theorem \ref{restr}  we have $L^{(2,1)}_r*L^{(2,1)}\subset L^{(2,\infty)}.$
On the other hand $L^1_r*L^1\subset L^1.$ By multilinear interpolation
theorem \cite[3.13.5, p. 76]{bl} we get
$L^{(p,s)}_r*L^{(p,t)}\subset L^{(p,u)}$ where $1\le p<2$ and
$1+1/u=1/s+1/t.$ Taking $u=\infty,$ $t=p$ and $s=p'$ gives
$L^{(p,p')}_r*L^p\subset L^{(p,\infty)}.$ This gives the first inequality.

On the other hand for $f=\sum_{n=0}^\infty f_n\chi_n$
by  (\ref{weak}) and by duality (\ref{dual})  we have
$$\|\lambda(f)\|_{p\to (p,\infty)}=\|\lambda(f)\|_{(p',1)\to p'}
\ge c\,
\sup_n\,q^{-n/p'}\|f*\chi_n\|_{p'}.
$$
Similarly as in the proof of Theorem \ref{weak} we obtain
$$f*\chi_n\ge \sum_{l=0}^\infty q^{(n-l)/2}\left[
\sum^{l+n}_{m=|n-l|\atop m\equiv l+n\mod 2}q^{m/2}f_m\right]\chi_l
$$
Hence
\begin{multline*}
q^{-n}\|f*\chi_n\|_{p'}^{p'}\ge \sum_{l=0}^nq^{p'(n-l)/2}q^{l-n}
\left[
\sum^{l+n}_{m=n-l\atop m\equiv l+n\mod 2}q^{m/2}f_m\right]^{p'}\\
\ge \sum_{l=0}^n q^{(n-l)(p'-1)}f_{n-l}^{p'}= \sum_{l=0}^n
q^{lp'/p}f_{l}^{p'}.
\end{multline*}
Taking supremum with respect to $n$ and raising to the power $1/p'$  give
$$\|\lambda(f)\|_{p\to (p,\infty)}\ge c\left (\sum_{n=0}^\infty f_n^{p'}q^{np'/p}\right
)^{1/p'}.$$  Since the norm of $f=\sum_{n=0}^\infty f_n\chi_n$ in
$L^{(p,p')}_r$ is equivalent to $ \left (\sum_{n=0}^\infty f_n^{p'}q^{np'/p}\right )^{1/p'}$
the second inequality is proved.
\qed
\section{Other estimates}
\begin{thm}
For $1\le s\le 2\le t\le \infty$ we have
$$cn^{1-1/s+1/t}q^{n/2}\le\|\lambda(\chi_n)\|_{(2,s)\to (2,t)}\le Cn^{1-1/s+1/t}q^{n/2}.$$
\end{thm}

{\it Proof.}
In order to get the second inequality we use only interpolation.
First observe that the inequality is valid for $s=2,\ t=\infty$ by Theorem \ref{weak} and for
$s=t=2$ by \cite{co, p0}. Hence by complex interpolation of the Lorentz spaces
it is valid for $s=2,\, t\ge 2.$

Next it is valid for $s=1,\, t=\infty$ by Theorem \ref{restr} and for $s=t=2.$
Hence by complex interpolation it is valid for $1\le s\le 2,\, t=s'.$

Now we can use again complex interpolation to get the conclusion for $1\le
s\le 2 \le t\le \infty. $

The estimate from below can be obtained from
$$\|\lambda(\chi_n)\|_{(2,s)\to (2,t)}\ge  {\|\chi_n*f\|_{(2,t)}\over
\|f\|_{(2,s)}},$$
where $f=\sum_{k=0}^{2n}q^{-k/2}\chi_k.$
\qed

Theorems 1, 2  and 5 suggest the following.

\noindent{\bf Conjecture.} Let $f=\sum_{n=0}^\infty f_n\chi_n\ge 0.$ Then for $1\le s\le
2$ the operator $\lambda(f)$ maps $L^{(2,s)} $ into $L^{(2,\infty)}$ if
and only if
$$ \sum_{n,m=0}^\infty
f_nf_mq^{-(n+m)/2}\{1+\min(n^{1/s'},m^{1/s'})\}<+\infty.$$

\end{document}